\UseRawInputEncoding 

\documentclass{article}

\usepackage{latexsym}
\usepackage{amsfonts}
\usepackage{amssymb}
\usepackage{cite}
\usepackage{color}

\newcommand{\supp}{{\rm supp}}
\newcommand{\mS}{\mathbb{S}}
\newcommand{\mK}{\mathbb{K}}

\newcommand{\brem}{\begin{remark}}
\newcommand{\erem}{\end{remark}}
\newcommand{\blem}{\begin{lemma}}
\newcommand{\elem}{\end{lemma}}
\newcommand{\bth}{\begin{theorem}}
\newcommand{\ethm}{\end{theorem}}
\newcommand{\benu}{\begin{enumerate}}
\newcommand{\eenu}{\end{enumerate}}
\newcommand{\bdes}{\begin{description}}
\newcommand{\edes}{\end{description}}
\newcommand{\bdf}{\begin{definition}}
\newcommand{\edf}{\end{definition}}
\newcommand{\bcor}{\begin{cor}}
\newcommand{\ecor}{\end{cor}}
\newcommand{\bprp}{\begin{proposition}}
\newcommand{\eprp}{\end{proposition}}
\newcommand{\bmlem}{\begin{mlemma}}
\newcommand{\emlem}{\end{mlemma}}
\newcommand{\bclm}{\begin{claim}}
\newcommand{\eclm}{\end{claim}}
\newcommand{\bprf}{{\bf Proof}.\hspace{2mm}}

\newcommand{\eprf}{\hspace*{\fill} $\Box$}

\newcommand{\beqn}{\begin{equation}}
\newcommand{\eeqn}{\end{equation}}
\newcommand{\beqnarr}{\begin{eqnarray}}
\newcommand{\eeqnarr}{\end{eqnarray}}
\newcommand{\beqnarrs}{\begin{eqnarray*}}
\newcommand{\eeqnarrs}{\end{eqnarray*}}

\newcommand{\spand}{\,\&\,}

\newtheorem{theorem}{Theorem}[section]
\newtheorem{definition}[theorem]{Definition}
\newtheorem{proposition}[theorem]{Proposition}
\newtheorem{lemma}[theorem]{Lemma}
\newtheorem{cor}[theorem]{Corollary}
\newtheorem{remark}[theorem]{Remark}
\newtheorem{mlemma}[theorem]{Main Lemma}
\newtheorem{claim}[theorem]{Claim}

\newcommand{\alp}{\alpha}

\newcommand{\del}{\delta}
\newcommand{\Del}{\Delta}
\newcommand{\ome}{\omega}
\newcommand{\Ome}{\Omega}
\newcommand{\bet}{\beta}
\newcommand{\gam}{\gamma}

\newcommand{\kap}{\kappa}
\newcommand{\sig}{\sigma}
\newcommand{\Sig}{\Sigma}
\newcommand{\tht}{\theta}

\newcommand{\Lam}{\Lambda}
\newcommand{\vphi}{\varphi}

\newcommand{\fal}{\forall}
\newcommand{\exi}{\exists}
\newcommand{\rarw }{\rightarrow}
\newcommand{\Rarw }{\Rightarrow}

\newcommand{\lrarw}{\leftrightarrow}
\newcommand{\Lrarw}{\Leftrightarrow}

\newcommand{\calc}{{\cal C}}

\newcommand{\calg}{{\cal G}}
\newcommand{\calh}{{\cal H}}

\newcommand{\calw}{{\cal W}}

\title{Well-foundedness proof for $\Pi^{1}_{1}$-reflection
}

\author{Toshiyasu Arai
\\
Graduate School of Mathematical Sciences,
University of Tokyo
\\
3-8-1 Komaba, Meguro-ku,
Tokyo 153-8914, JAPAN
\\
tosarai@ms.u-tokyo.ac.jp
}

\date{}

\begin{document}

\maketitle

\begin{abstract}
In \cite{ghent} it is shown that an ordinal
$\psi_{\Omega}(\varepsilon_{\mathbb{S}^{+}+1})$ is an upper bound for the proof-theoretic ordinal of a set theory ${\sf KP}\ome+(M\prec_{\Sigma_{1}}V)$.
In this note we show that 
${\sf KP}\ome+(M\prec_{\Sigma_{1}}V)$
 proves the well-foundedness up to $\psi_{\Omega}(\ome_{n}(\mathbb{S}^{+}+1))$ for each $n$.
\end{abstract}

\maketitle

\section{Introduction}\label{sect:introduction}
In \cite{ghent} the following theorem is shown, where
${\sf KP}\ome+(M\prec_{\Sigma_{1}}V)$ extends  ${\sf KP}\ome$ 
with an axiom stating that
`there exists an non-empty and transitive set $M$ such that $M\prec_{\Sigma_{1}}V$'.
$\Omega=\omega_{1}^{CK}$ and $\psi_{\Omega}$ is a collapsing function such that $\psi_{\Omega}(\alpha)<\Omega$.
$\mathbb{S}$ is an ordinal term denoting a stable ordinal, and 
$\mathbb{S}^{+}$ the least admissible ordinal above $\mathbb{S}$ in the theorems.

\begin{theorem}\label{thm:2}
Suppose ${\sf KP}\ome+(M\prec_{\Sigma_{1}}V)\vdash\theta^{L_{\Omega}}$
for a $\Sigma_{1}$-sentence $\theta$.
Then we can find an $n<\omega$ such that for $\alpha=\psi_{\Omega}(\omega_{n}(\mathbb{S}^{+}+1))$,
$L_{\alpha}\models\theta$.
\end{theorem}

$OT$ denotes a computable notation system of ordinals in \cite{singlewfprf} for
an ordinal analysis of ${\sf KP}\ell^{r}+(M\prec_{\Sigma_{1}}V)$, or equivalently 
of $\Sigma^{1-}_{2}\mbox{-CA}+\Pi^{1}_{1}\mbox{-CA}_{0}$.
$OT_{N}$ is a restriction of $OT$ such that $OT=\bigcup_{0<N<\omega}OT_{N}$ and
$\psi_{\Omega}(\varepsilon_{\Omega_{\mathbb{S}+N}+1})$ denotes
the order type of $OT_{N}\cap\Omega$.
Let $OT(\Pi^{1}_{1})=OT_{1}$.
The aim of this paper is to show the following theorem, thereby the bound in Theorem \ref{thm:2}
is seen to be tight.

\begin{theorem}\label{th:wf}
${\sf KP}\ome+(M\prec_{\Sigma_{1}}V)$
proves the well-foundedness up to $\psi_{\Omega}(\ome_{n}(\mathbb{S}^{+}+1))$ for 
{\rm each} $n$.
\end{theorem}

The ordinal $\psi_{\Omega}(\varepsilon_{\mathbb{S}^{+}+1})$ 
is the proof-theoretic ordinal of 
${\sf KP}\ome+(M\prec_{\Sigma_{1}}V)$.

\begin{theorem}\label{th:main}
$\psi_{\Omega}(\varepsilon_{\mathbb{S}^{+}+1})
= |{\sf KP}\ome+(M\prec_{\Sigma_{1}}V)|_{\Sigma_{1}^{\Omega}}.
$
\end{theorem}

To prove the well-foundedness of a computable notation system,
we utilize
the distinguished class introduced by W. Buchholz\cite{Buchholz75}.

A set theory ${\sf KP}\ome+(M\prec_{\Sigma_{1}}V)$ extends ${\sf KP}\ome$
 by adding an individual constant $M$ and
the axioms for the constant $M$:
$M$ is non-empty $M\neq\emptyset$, transitive $\forall x\in M\forall y\in x(y\in M)$, and stable $M\prec_{\Sigma_{1}}V$ for the universe $V$. $M\prec_{\Sigma_{1}}V$ means that
$
 \varphi(u_{1},\ldots,u_{n}) \land \{u_{1},\ldots,u_{n}\}\subset M \to \varphi^{M}(u_{1},\ldots,u_{n})
$
for each $\Sigma_{1}$-formula $\varphi$ in the set-theoretic language.

Since the axiom $\bet$ does not hold in the theory ${\sf KP}\ome+(M\prec_{\Sigma_{1}}V)$,
we need to modify the proof in \cite{singlewfprf}, cf.\,subsection \ref{subsec:disting}.
Proofs of propositions and lemmas are omitted when they are found in \cite{Wienpi3d,singlewfprf}.

\section{Ordinals for one stable ordinal}\label{sect:ordinalnotation}
In this section let us recall briefly ordinal notations systems in \cite{singlewfprf}.

For ordinals $\alpha\geq\beta$,
$\alpha-\beta$ denotes the ordinal $\gamma$ such that $\alpha=\beta+\gamma$.
Let $\alpha$ and $\beta$ be ordinals.
$\alpha\dot{+}\beta$ denotes the sum $\alpha+\beta$
when $\alpha+\beta$ equals to the commutative (natural) sum $\alpha\#\beta$, i.e., when
either $\alpha=0$ or $\alpha=\alpha_{0}+\omega^{\alpha_{1}}$ with
$\omega^{\alpha_{1}+1}>\beta$.

$\mathbb{S}$ denotes a weakly inaccessible cardinal, and
$\Lam=\mathbb{S}^{+}$ the next regular cardinal above $\mathbb{S}$.

\bdf\label{df:Lam}
{\rm
Let $\Lambda=\mathbb{S}^{+}$.
$\varphi_{b}(\xi)$ denotes the binary Veblen function on 
$\Lambda^{+}$ with $\varphi_{0}(\xi)=\omega^{\xi}$, and
$\tilde{\varphi}_{b}(\xi):=\varphi_{b}(\Lambda\cdot \xi)$ for the epsilon number 
$\Lambda$.

Let $b,\xi<\Lambda^{+}$.
$\theta_{b}(\xi)$ [$\tilde{\theta}_{b}(\xi)$] denotes
a $b$-th iterate of $\varphi_{0}(\xi)=\omega^{\xi}$ [of $\tilde{\varphi}_{0}(\xi)=\Lambda^{\xi}$], resp.
}
\edf

\bdf\label{df:Lam2}
{\rm
Let $\xi<\varphi_{\Lambda}(0)$ be a non-zero ordinal with its normal form:
\begin{equation}\label{eq:CantornfLam}
\xi=\sum_{i\leq m}\tilde{\theta}_{b_{i}}(\xi_{i})\cdot a_{i}=_{NF}
\tilde{\theta}_{b_{m}}(\xi_{m})\cdot a_{m}+\cdots+\tilde{\theta}_{b_{0}}(\xi_{0})\cdot a_{0}
\end{equation}
where
$\tilde{\theta}_{b_{i}}(\xi_{i})>\xi_{i}$,
$\tilde{\theta}_{b_{m}}(\xi_{m})>\cdots>\tilde{\theta}_{b_{0}}(\xi_{0})$, 
$b_{i}=\omega^{c_{i}}<\Lambda$, and
$0<a_{0},\ldots,a_{m}<\Lambda$.
$SC_{\Lam}(\xi)=\bigcup_{i\leq m}(\{a_{i}\}\cup SC_{\Lam}(\xi_{i}))$.

$\tilde{\theta}_{b_{0}}(\xi_{0})$ is said to be the \textit{tail} of $\xi$, denoted 
$\tilde{\theta}_{b_{0}}(\xi_{0})=tl(\xi)$, and
$\tilde{\theta}_{b_{m}}(\xi_{m})$ the \textit{head} of $\xi$, denoted 
$\tilde{\theta}_{b_{m}}(\xi_{m})=hd(\xi)$.

\begin{enumerate}
\item\label{df:Exp2.3}
 $\zeta$ is a \textit{segment} of $\xi$
 iff there exists an $n\, (0\leq n\leq m+1)$
 such that
 $\zeta=_{NF}\sum_{i\geq n}\tilde{\theta}_{b_{i}}(\xi_{i})\cdot a_{i}=
 \tilde{\theta}_{b_{m}}(\xi_{m})\cdot a_{m}+\cdots+\tilde{\theta}_{b_{n}}(\xi_{n})\cdot a_{n}$
 for $\xi$ in (\ref{eq:CantornfLam}).

\item\label{df:thtminus}
Let $\zeta=_{NF}\tilde{\theta}_{b}(\xi)$ with $\tilde{\theta}_{b}(\xi)>\xi$ and $b=\omega^{b_{0}}$,
and $c$ be ordinals.
An ordinal $\tilde{\theta}_{-c}(\zeta)$ is defined recursively as follows.
If $b\geq c$, then $\tilde{\theta}_{-c}(\zeta)=\tilde{\theta}_{b-c}(\xi)$.
Let $c>b$.
If $\xi>0$, then
$\tilde{\theta}_{-c}(\zeta)=\tilde{\theta}_{-(c-b)}(\tilde{\theta}_{b_{m}}(\xi_{m}))$ for the head term 
$hd(\xi)=\tilde{\theta}_{b_{m}}(\xi_{m})$ of 
$\xi$ in (\ref{eq:CantornfLam}).
If $\xi=0$, then let $\tilde{\theta}_{-c}(\zeta)=0$.

\end{enumerate}
}
\edf

\bdf\label{df:Lam3}
{\rm

  \begin{enumerate}
  \item
A function $f:\Lam\to \varphi_{\Lam}(0)$ with a \textit{finite} support
${\rm supp}(f)=\{c<\Lam: f(c)\neq 0\}\subset \Lam$ is said to be a \textit{finite function}
if
$\forall i>0(a_{i}=1)$ and $a_{0}=1$ when $b_{0}>1$
in
$f(c)=_{NF}\tilde{\theta}_{b_{m}}(\xi_{m})\cdot a_{m}+\cdots+\tilde{\theta}_{b_{0}}(\xi_{0})\cdot a_{0}$
for any $c\in{\rm supp}(f)$.

It is identified with the finite function $f\!\upharpoonright\! {\rm supp}(f)$.
When $c\not\in {\rm supp}(f)$, let $f(c):=0$.
$SC_{\Lam}(f):=\bigcup\{\{c\}\cup SC_{\Lam}(f(c))\}: c\in {\rm supp}(f)\}$.
$f,g,h,\ldots$ range over finite functions.

For an ordinal $c$, $f_{c}$ and $f^{c}$ are restrictions of $f$ to the domains
${\rm supp}(f_{c})=\{d\in{\rm supp}(f): d< c\}$ and ${\rm supp}(f^{c})=\{d\in{\rm supp}(f): d\geq c\}$.
$g_{c}*f^{c}$ denotes the concatenated function such that
${\rm supp}(g_{c}*f^{c})={\rm supp}(g_{c})\cup {\rm supp}(f^{c})$, 
$(g_{c}*f^{c})(a)=g(a)$ for $a<c$, and
$(g_{c}*f^{c})(a)=f(a)$ for $a\geq c$.


\item\label{df:Exp2.5}
Let $f$ be a finite function and $c,\xi$ ordinals.
A relation $f<^{c}\xi$ is defined by induction on the
cardinality of the finite set $\{d\in {\rm supp}(f): d>c\}$ as follows.
If $f^{c}=\emptyset$, then $f<^{c}\xi$ holds.
For $f^{c}\neq\emptyset$,
$f<^{c}\xi$ iff
there exists a segment $\mu$ of $\xi$ such that
$f(c)< \mu$
and 
$f<^{c+d} \tilde{\theta}_{-d}(tl(\mu))$ 
for $d=\min\{c+d\in {\rm supp}(f): d>0\}$.


\end{enumerate}

}
\edf



\bprp\label{prp:idless}
$f<^{c}\xi\leq\zeta \Rightarrow f<^{c}\zeta$.
\eprp

In the following Definition \ref{df:Cpsiregularsm}, 
$\varphi\alpha\beta=\varphi_{\alpha}(\beta)$ denotes the binary Veblen function on $\Lam^{+}=\mathbb{S}^{++}$,
$\tilde{\theta}_{b}(\xi)$ the function defined in Definition \ref{df:Lam}
for $\Lambda=\mathbb{S}^{+}$.
For $\alp<\mS$, $\alp^{+}$ denotes the next regular cardinal above $\alp$.

For 
$a<\varepsilon_{\Lam+1}$,
$c<\Lam$, and
$\xi<\Gamma_{\Lam+1}$, 
define simultaneously 
classes $\mathcal{H}_{a}(X)\subset\Gamma_{\Lam+1}$,
$Mh^{a}_{c}(\xi)\subset(\mathbb{S}+1)$, and 
ordinals $\psi_{\kappa}^{f}(a)\leq\kappa$ by recursion on ordinals $a$ as follows.

\begin{definition}\label{df:Cpsiregularsm}
{\rm
Let
$\Lam=\mathbb{S}^{+}$.
Let $a<\varepsilon_{\Lam+1}$ and $X\subset\Gamma_{\Lam+1}$.

\begin{enumerate}
\item\label{df:Cpsiregularsm.1}
(Inductive definition of $\mathcal{H}_{a}(X)$.)

\begin{enumerate}
\item\label{df:Cpsiregularsm.10}
$\{0,\Omega_{1},\mathbb{S},\mS^{+}\}\cup X\subset\mathcal{H}_{a}(X)$.

\item\label{df:Cpsiregularsm.11}
If $x, y \in \mathcal{H}_{a}(X)$,
then $x+y\in \mathcal{H}_{a}(X)$,
and 
$\varphi xy\in \mathcal{H}_{a}(X)$.

\item\label{df:Cpsiregularsm.12}
Let $\alpha\in\mathcal{H}_{a}(X)\cap\mathbb{S}$. Then 
$\alp^{+}\in\mathcal{H}_{a}(X)$.

\item\label{df:Cpsiregularsm.1345}
Let $\alpha=\psi_{\pi}^{f}(b)$ with $\{\pi,b\}\subset\mathcal{H}_{a}(X)$, 
$b<a$, and a finite function $f$ such that
$SC_{\Lam}(f)\subset\mathcal{H}_{a}(X)\cap\mathcal{H}_{b}(\alpha)$.
Then $\alpha\in\mathcal{H}_{a}(X)$.

\end{enumerate}

\item\label{df:Cpsiregularsm.2}
 (Definitions of $Mh^{a}_{c}(\xi)$ and $Mh^{a}_{c}(f)$)
\\
The classes $Mh^{a}_{c}(\xi)$ are defined for $c< \Lam$,
and ordinals $a<\varepsilon_{\Lam+1}$, $\xi<\Gamma_{\Lam+1}$.
Let $\pi$ be a regular ordinal$\leq \mathbb{S}$. Then 
by main induction on ordinals $\pi\leq\mathbb{S}$
with subsidiary induction on $c<\Lam$ 
we define $\pi\in Mh^{a}_{c}(\xi)$ iff 
$\{a,c,\xi\}\subset\mathcal{H}_{a}(\pi)$ and
\begin{equation}\label{eq:dfMhkh}
 \forall f<^{c}\xi 
 \forall g \left(
SC_{\Lam}(f,g) \subset\mathcal{H}_{a}(\pi) 
 \,\&\, 
\pi\in Mh^{a}_{0}(g_{c})
 \Rightarrow \pi\in M(Mh^{a}_{0}(g_{c}*f^{c}))
 \right)
\end{equation}
where $f, g$ vary through 
finite
 functions,
and 
\begin{eqnarray*}
Mh^{a}_{c}(f)  & := & \bigcap\{Mh^{a}_{d}(f(d)): d\in {\rm supp}(f^{c})\}
\\
& = &
\bigcap\{Mh^{a}_{d}(f(d)): c\leq d\in {\rm supp}(f)\}.
\end{eqnarray*}
In particular
$Mh^{a}_{0}(g_{c})=\bigcap\{Mh^{a}_{d}(g(d)): d\in {\rm supp}(g_{c})\}
=\bigcap\{Mh^{a}_{d}(g(d)): c> d\in {\rm supp}(g)\}$.
When $f=\emptyset$ or $f^{c}=\emptyset$, let $Mh^{a}_{c}(\emptyset):=\Lam$.

\item\label{df:Cpsiregularsm.3}
 (Definition of $\psi_{\pi}^{f}(a)$)
\\
 Let $a<\varepsilon_{\Lam+1}$ 
 be an ordinal, $\pi$ a regular ordinal and
 $f$ a finite function.
Then let
\begin{equation}\label{eq:Psivec}
\hspace{-10mm}
\psi_{\pi}^{f}(a)
 :=  \min(\{\pi\}\cup\{\kappa\in Mh^{a}_{0}(f)\cap\pi:   \mathcal{H}_{a}(\kappa)\cap\pi\subset\kappa ,
   \{\pi,a\}\cup SC_{\Lam}(f)\subset\mathcal{H}_{a}(\kappa)
\})
\end{equation}
For the empty function $\emptyset$,
$\psi_{\pi}(a):=\psi_{\pi}^{\emptyset}(a)$.

\item
For classes $A\subset(\mathbb{S}+1)$, let
$\alpha\in M^{a}_{c}(A)$ iff $\alpha\in A$ and
\begin{equation}\label{eq:Mca}
\forall g
[
\alpha\in Mh_{0}^{a}(g_{c}) \,\&\, SC_{\Lam}(g_{c})\subset\mathcal{H}_{a}(\alpha) \Rightarrow
\alpha\in M\left( Mh_{0}^{a}(g_{c}) \cap A \right)
]
\end{equation}
\end{enumerate}
}

\end{definition}

Assuming an existence of a 
shrewd cardinal introduced by M. Rathjen\cite{RathjenAFML2},
we show in \cite{singlestable} that
$\psi_{\mathbb{S}}^{f}(a)<\mathbb{S}$
if $\{a,c,\xi\}\subset\mathcal{H}_{a}(\mathbb{S})$ with 
$c<\mathbb{S}^{+}$, $a,\xi<\varepsilon_{\mathbb{S}^{+}+1}$,
and ${\rm supp}(f)=\{c\}$ and $f(c)=\xi$.
Moreover $\psi_{\pi}^{g}(b)<\pi$ 
provided that
$\pi\in Mh^{b}_{0}(f)$, $SC_{\Lam}(g)\cup\{\pi,b\}\subset\mathcal{H}_{b}(\pi)$, and
$g$ is a finite function defined from a finite function $f$ and ordinals $d,c$ as follows.
$d<c\in \supp(f)$ 
with
$(d,c)\cap \supp(f)=(d,c)\cap \supp(g)=\emptyset$, 
$g_{d}=f_{d}$,
$g(d)<f(d)+\tilde{\theta}_{c-d}(f(c))\cdot\omega$, and
$g<^{c}f(c)$.
Also the following Lemma \ref{lem:stepdown} is shown in \cite{singlestable}.

\begin{lemma}\label{lem:stepdown}
Assume $\mathbb{S}\geq\pi\in Mh^{a}_{d}(\xi)\cap Mh^{a}_{c}(\xi_{0})$, $\xi_{0}\neq 0$,
$d<c$, and
$\{a,c,d\}\subset\mathcal{H}_{a}(\pi)$. 
Moreover let $\tilde{\theta}_{c-d}(\xi_{0})\geq\xi_{1}\in\mathcal{H}_{a}(\pi)$
and $tl(\xi)>\xi_{1}$ when $\xi\neq 0$.
Then
$\pi\in Mh^{a}_{d}(\xi+\xi_{1})\cap M^{a}_{d}(Mh^{a}_{d}(\xi+\xi_{1}))$.
\end{lemma}

\subsection{Normal forms in ordinal notations}

\begin{definition}\label{df:irreducible}
{\rm

An \textit{irreducibility} of finite functions $f$ is defined by induction on the cardinality
$n$ of the finite set $\supp(f)$.
If $n\leq 1$, $f$ is defined to be irreducible.
Let $n\geq 2$ and $c<c+d$ be the largest two elements in $\supp(f)$, and let $g$ be 
a finite function
such that $\supp(g)=\supp(f_{c})\cup\{c\}$, $g_{c}=f_{c}$ and
$g(c)=f(c)+\tilde{\theta}_{d}(f(c+d))$.
Then $f$ is irreducible iff 
$tl(f(c))> \tilde{\theta}_{d}(f(c+d))$ and
$g$ is irreducible.

}
\end{definition}

\begin{definition}\label{df:lx}
 {\rm 
 Let  $f,g$ be irreducible functions, and $b,a$ ordinals.
 \begin{enumerate}
 \item\label{df:lxx}
Let us define a relation $f<^{b}_{lx}g$
by induction on the cardinality of the finite set
$\{e\in{\rm supp}(f)\cup{\rm supp}(g): e\geq b\}$ as follows.
$f<^{b}_{lx}g$ holds iff $f^{b}\neq g^{b}$ and
for the ordinal $c=\min\{c\geq b : f(c)\neq g(c)\}$,
one of the following conditions is met:

\begin{enumerate}

\item\label{df:lx.23}
$f(c)<g(c)$ and let $\mu$ be the shortest segment of $g(c)$ such that $f(c)<\mu$.
Then for any $c<c+d\in{\rm supp}(f)$,  
if $tl(\mu)\leq\tilde{\theta}_{d}(f(c+d))$, then 
$f<_{lx}^{c+d}g$ holds.

\item\label{df:lx.24}
$f(c)>g(c)$ and let $\nu$ be the shortest segment of $f(c)$ such that $\nu>g(c)$.
Then there exist a $c<c+d\in {\rm supp}(g)$ such that
$f<_{lx}^{c+d}g$ and
$tl(\nu)\leq \tilde{\theta}_{d}(g(c+d))$.

\end{enumerate}


\item
$Mh^{a}_{b}(f)\prec Mh^{a}_{b}(g)$ holds iff
\[
\forall\pi\in Mh^{a}_{b}(g)\forall b_{0}\leq b
\left(
SC_{\Lam}(f)\subset\mathcal{H}_{a}(\pi) \,\&\, \pi\in Mh^{a}_{b_{0}}(f_{b})
 \Rightarrow \pi\in M(Mh^{a}_{b_{0}}(f))
\right)
.
\]
\end{enumerate}
}
\end{definition}


\begin{lemma}\label{lem:psinucomparison}
Let $f,g$ be irreducible finite functions, and $b$ an ordinal such that $f^{b}\neq g^{b}$.
If $f<^{b}_{lx}g$, then
$Mh^{a}_{b}(f)\prec Mh^{a}_{b}(g)$ holds for every ordinal $a$.
\end{lemma}

\begin{proposition}\label{prp:psicomparison}
Let $f,g$ be 
irreducible finite functions, and assume that
$\psi_{\pi}^{f}(b)<\pi$ and $\psi_{\kappa}^{g}(a)<\kappa$.

Then $\psi_{\pi}^{f}(b)<\psi_{\kappa}^{g}(a)$ iff one of the following cases holds:
\benu
\item\label{prp:psicomparison.0}
$\pi\leq \psi_{\kappa}^{g}(a)$.

\item\label{prp:psicomparison.1}
$b<a$, $\psi_{\pi}^{f}(b)<\kappa$ and 
$SC_{\Lam}(f)\cup\{\pi,b\}\subset\mathcal{H}_{a}(\psi_{\kappa}^{g}(a))$.

\item\label{prp:psicomparison.2}
$b>a$ and $SC_{\Lam}(g)\cup\{\kappa,a\}\not\subset\mathcal{H}_{b}(\psi_{\pi}^{f}(b))$.

\item\label{prp:psicomparison.25}
$b=a$, $\kappa<\pi$ and $\kappa\not\in\mathcal{H}_{b}(\psi_{\pi}^{f}(b))$.

\item\label{prp:psicomparison.3}
$b=a$, $\pi=\kappa$, $SC_{\Lam}(f)\subset\mathcal{H}_{a}(\psi_{\kappa}^{g}(a))$, and
$f<^{0}_{lx}g$.

\item\label{prp:psicomparison.4}
$b=a$, $\pi=\kappa$, 
$SC_{\Lam}(g)\not\subset\mathcal{H}_{b}(\psi_{\pi}^{f}(b))$.

\eenu

\end{proposition}

\begin{definition}\label{df:Lam.of}
{\rm
\benu
\item
$a(\xi)$ denotes an ordinal defined recursively by
$a(0)=0$, and
$a(\xi)=\sum_{i\leq m}\tilde{\theta}_{b_{i}}(\omega\cdot a(\xi_{i}))$
when $\xi=_{NF}\sum_{i\leq m}\tilde{\theta}_{b_{i}}(\xi_{i})\cdot a_{i}$ in (\ref{eq:CantornfLam}).

\item
For irreducible functions $f$ let us associate ordinals $o(f)<\Gamma_{\mS^{+}+1}$ as follows.
$o(\emptyset)=0$ for the empty function $f=\emptyset$.
Let $\{0\}\cup \supp(f)=\{0=c_{0}<c_{1}<\cdots<c_{n}\}$,
$f(c_{i})=\xi_{i}<\Gamma_{\mS^{+}+1}$ for $i>0$, and $\xi_{0}=0$.
Define ordinals $\zeta_{i}=o(f;c_{i})$ by
$\zeta_{n}=\omega\cdot a(\xi_{n})$, and $\zeta_{i}=\omega\cdot a(\xi_{i})+\tilde{\theta}_{c_{i+1}-c_{i}}(\zeta_{i+1}+1)$.
Finally let $o(f)=\zeta_{0}=o(f;c_{0})$.

\item
Let $SC_{\Lam}(f)<\mu<\Lam$ be an epsilon number.
Then $o_{\mu}(f)$ is defined from $o(f)$ by replacing the base $\Lam$ of $\tilde{\tht}$
in $f(c)$ by $\mu$. This means that $\Lam$ is replaced by $\mu$, and
$\tilde{\tht}_{1}(\xi)=\Lam^{\xi}$ by $\mu^{\xi}$.

\eenu
}
\end{definition}

\begin{lemma}\label{lem:of}
Let $f$ be an irreducible finite function defined from an irreducible function $g$ and ordinals $c,d$
as follows.
$f_{c}=g_{c}$, 
$c<d\in \supp(g)$ with
$(c,d)\cap \supp(g)=
(c,d)\cap \supp(f)=
\emptyset$, 
$f(c)<g(c)+\tilde{\theta}_{d-c}(g(d))\cdot\omega$, and
$f<^{d}g(d)$.
Then $o(f)<o(g)$ holds.

Moreover when $SC_{\Lam}(f,g)<\mu<\Lam$,
$o_{\mu}(f)<o_{\mu}(g)$ holds.
\end{lemma}

\begin{lemma}\label{lem:oflx}
For irreducible finite functions $f$ and $g$, assume
$f<_{lx}^{0}g$.
Then $o(f)<o(g)$ holds.

Moreover when $SC_{\Lam}(f,g)<\mu<\Lam$,
$o_{\mu}(f)<o_{\mu}(g)$ holds.
\end{lemma}

By Proposition \ref{prp:psicomparison} a notation system $OT(\Pi^{1}_{1})=OT_{1}$
is defined.
\bdf\label{df:notationsystempi11}
{\rm
$OT(\Pi^{1}_{1})$ is closed under $\mS>\alp\mapsto\alp^{+}$.
There are two cases when an ordinal term $\psi_{\pi}^{f}(a)$ is constructed in $OT(\Pi^{1}_{1})$, 
from $\{\pi,a\}\subset OT(\Pi^{1}_{1})$ and an irreducible function $f$
with $SC_{\Lam}(f)\subset OT(\Pi^{1}_{1})$ and $\Lam=\mS^{+}$.
$E_{\mathbb{S}}(\alp)$ denotes the set of subterms$<\mS$ of $\alp$.
\benu
 \item\label{df:notationsystem.10}
Let $\xi,a,c\in OT(\Pi^{1}_{1})$,
$\xi>0$, $c<\mathbb{S}^{+}$ and
$\{\xi,a,c\}\subset\calh_{a}(\alp)$.
Then
$\alpha=\psi_{\mathbb{S}}^{f}(a)\in OT(\Pi^{1}_{1})$ 
and $\alp^{+}\in OT(\Pi^{1}_{1})$ 
with ${\rm supp}(f)=\{c\}$ and $f(c)=\xi$
if $\max(SC_{\mS^{+}}(f))\leq \max(SC_{\mS^{+}}(a))$.
Let $f=m(\alp)$.

 \item\label{df:notationsystem.11}
Let $\{a,d,\pi\}\subset OT(\Pi^{1}_{1})$, 
$f=m(\pi)$,
 $d<c\in \supp(f)$,
and $(d,c)\cap \supp(f)=\emptyset$.
Let $g$ be an irreducible function such that 
$SC_{\Lam}(g)=\bigcup\{\{c,g(c)\}: c\in {\rm supp}(g)\}\subset OT(\Pi^{1}_{1})$,
$g_{d}=f_{d}$, $(d,c)\cap \supp(g)=\emptyset$
$g(d)<f(d)+\tilde{\theta}_{c-d}(f(c))\cdot\omega$, 
and $g<^{c}f(c)$.
Moreover if $\max(SC_{\Lam}(f))<\mu<\Lam$ for an epsilon number $\mu$, then 
$\max(SC_{\Lam}(g))<\mu$.

Then 
$\alpha=\psi_{\pi}^{g}(a)\in OT(\Pi^{1}_{1})$ and $\alp^{+}\in OT(\Pi^{1}_{1})$ if 
$\{\pi,a\}\cup SC_{\Lam}(f,g)\subset\calh_{a}(\alp)$,
 and, cf.\,Proposition \ref{prp:psimS}.

\begin{equation}\label{eq:notationsystem.11}
SC_{\Lam}(g)\subset M_{\alp}
\end{equation}

\eenu
}
\edf

$M_{\alp}$ is defined as follows.

\bdf
{\rm
For ordinal terms $\psi_{\sig}^{f}(a)\in \Psi_{\mS}\subset OT(\Pi^{1}_{1})$,
define
$m(\psi_{\sig}^{f}(a)):=f$ and 
${\tt p}_{0}(\psi_{\sig}^{f}(a))={\tt p}_{0}(\sig)$ if $\sig<\mS$, and
${\tt p}_{0}(\psi_{\mS}^{f}(a))=a$.
}
\edf

\bdf
$M_{\rho}:=\calh_{b}(\rho)$ for $b={\tt p}_{0}(\rho)$ and $\rho\in\Psi_{\mS}$.
\edf

\bdf\label{df:Lammu}
{\rm
For $\gam\prec\mS$, an epsilon number $\mS<\mu=\Lam(\gam)<\mS^{+}$ is defined.
Let $\gam=\psi_{\sig}^{f}(\alp)\preceq\psi_{\mS}^{g}(b)$ with $b={\tt p}_{0}(\gam)$.
Then $\Lam(\gam)$ denotes the least epsilon number $\mS<\mu<\mS^{+}$ such that
$\max(SC_{\mS^{+}}(b))<\mu$.
}
\edf
From Definition \ref{df:notationsystempi11} we see $\max(SC_{\mS^{+}}(f))<\Lam(\gam)$.

$OT(\Pi^{1}_{1})$ is closed under
$\alpha\mapsto\alpha[\rho/\mS]$ for $\alp\in M_{\rho}$.
Specifically if $\{\alp,\rho\}\subset OT(\Pi^{1}_{1})$ with
$\alp\in M_{\rho}$ and $\rho\in\Psi_{\mS}$, then
$\alp[\rho/\mS]\in OT(\Pi^{1}_{1})$.

\bdf\label{df:Mostwskicollaps}
{\rm
Let $\alpha\in M_{\rho}$ with $\rho\in \Psi_{\mS}$.
We define an ordinal $\alpha[\rho/\mS]$ recursively as follows.
$\alpha[\rho/\mS]:=\alpha$ when $\alpha<\mathbb{S}$.
In what follows assume $\alpha\geq\mathbb{S}$.

$\mathbb{S}[\rho/\mS]:=\rho$.
$\mK[\rho/\mS]\equiv(\mathbb{S}^{+})[\rho/\mS]:=\rho^{+}$.
$\left(\psi_{\mK}(a)\right)[\rho/\mS]=\left(\psi_{\mathbb{S}^{+}}(a)\right)[\rho/\mS]=\psi_{\rho^{+}}(a[\rho/\mS])$.
The map commutes with $+$ and $\varphi$.
}
\edf


\blem\label{lem:Mostowskicollapspi11}
For $\rho\in\Psi_{\mS}$,
$\{\alpha[\rho/\mS]:\alpha\in M_{\rho}\}$ is a transitive collapse of $M_{\rho}$ in the sense that
$\beta<\alpha\Leftrightarrow\beta[\rho/\mS]<\alpha[\rho/\mS]$,
$\beta\in\mathcal{H}_{\alpha}(\gamma)\Leftrightarrow 
\beta[\rho/\mS]\in\mathcal{H}_{\alpha[\rho/\mS]}(\gamma[\rho/\mS]))
$ for $\gamma>\mathbb{S}$,
and
$OT(\Pi^{1}_{1})\cap\alpha[\rho/\mS]=\{\beta[\rho/\mS]:\beta\in M_{\rho}\cap\alpha\}$
for $\alpha,\beta,\gamma\in M_{\rho}$.
\elem

\bprp\label{prp:EK2pi11}
Let $\rho\in\Psi_{\mS}$.
\benu
\item\label{prp:EK2.1}
$\mathcal{H}_{\gamma}(M_{\rho})\subset M_{\rho}$ 
if $\gam\leq{\tt p}_{0}(\rho)$.

\item\label{prp:EK2.2}
$M_{\rho}\cap\mS=\rho$ and $\rho\not\in M_{\rho}$.

\item\label{prp:EK2.3}
If $\sig<\rho$ and ${\tt p}_{0}(\sig)\leq {\tt p}_{0}(\rho)$, then
$M_{\sig}\subset M_{\rho}$.
\eenu
\eprp

\section{Well-foundedness proof with the maximal distinguished set}\label{sec:distinguished}
In this section working in the set theory ${\sf KP}\ome+(M\prec_{\Sigma_{1}}V)$, 
we show the well-foundedness of the notation system $OT(\Pi_{1}^{1})$ up to each
$\psi_{\Ome}(\ome_{n}(\mS^{+}+1))$.
Let us write $L_{\mS}$ for $M$, i.e., $L_{\mS}\prec_{\Sig_{1}}L$.
The proof is based on distinguished classes, which was first
introduced by Buchholz\cite{Buchholz75}.

\subsection{Distinguished sets}\label{subsec:disting}

 $X,Y,\ldots$ range over subsets of $OT(\Pi_{1}^{1})$. 
We define sets $\mathcal{C}^{\alpha}(X)\subset OT(\Pi_{1}^{1})$ for $\alpha\in OT(\Pi_{1}^{1})$ and
$X\subset OT(\Pi_{1}^{1})$ as follows.

\begin{definition}\label{df:CX}
{\rm 
Let $\alpha,\beta\in OT(\Pi_{1}^{1})$ and $X\subset OT(\Pi_{1}^{1})$.

$\mathcal{C}^{\alpha}(X)$ denotes the closure of $\{0,\Omega,\mathbb{S},\mathbb{S}^{+} \}\cup(X\cap\alpha)$
under $+$, 
$\sigma\mapsto \sig^{+}$,
$(\beta,\gamma)\mapsto\varphi\beta\gamma$,
and $(\sigma,\beta,f)\mapsto \psi_{\sigma}^{f}(\beta)$
for $ \sigma>\alpha$ in  $OT(\Pi_{1}^{1})$.

The last clause says that,
$\psi_{\sigma}^{f}(\beta)\in\mathcal{C}^{\alpha}(X)$ if
$\{\sigma,\beta\}\cup SC_{\Lam}(f)\subset\mathcal{C}^{\alpha}(X)$ and $\sigma>\alpha$.

}
\end{definition}

\begin{proposition}\label{lem:CX2} 
Assume $\forall\gamma\in X[ \gamma\in\mathcal{C}^{\gamma}(X)]$ for a set $X\subset OT(\Pi_{1}^{1})$.

\benu
\item\label{lem:CX2.3} 
$\alpha\leq\beta \Rarw \mathcal{C}^{\beta}(X)\subset \mathcal{C}^{\alpha}(X)$.

\item\label{lem:CX2.4} 
$\alpha<\beta<\alpha^{+} \Rarw \mathcal{C}^{\beta}(X)=\mathcal{C}^{\alpha}(X)$.
\eenu
\end{proposition}

\begin{definition}\label{df:wftg}
{\rm
\benu
\item 
$Prg[X,Y] :\Lrarw \forall\alpha\in X(X\cap\alpha\subset Y \to \alpha\in Y)$.

\item 
For a definable class $\mathcal{X}$, $TI[\mathcal{X}]$ denotes the schema:\\
$TI[\mathcal{X}] :\Lrarw Prg[\mathcal{X},\mathcal{Y}]\to \mathcal{X}\subset\mathcal{Y} \mbox{ {\rm holds for} any definable classes } \mathcal{Y}$.

\item
For $X\subset OT(\Pi_{1}^{1})$, $W(X)$ denotes the \textit{well-founded part} of $X$. 

\item 
$Wo[X] : \Lrarw X\subset W(X)$.

\item 
$\alp\in W_{\Sig}(X)$ denotes a $\Sig_{1}$-formula saying that
$\alp\in X$ and 
 `there exists an embedding
$f:X\cap(\alp+1)\to ON$', i.e.,
$\exi f\in{}^{\ome}ON\fal \bet,\gam\in X\cap(\alp+1)(\bet<\gam \to f(\bet)<f(\gam))$,
where $ON$ is the class of all ordinals, $\bet<\gam$ in $OT(\Pi^{1}_{1})$
and $f(\bet)<f(\gam)$ in $ON$.

\item 
$Wo_{\Sig}[X]$ denotes a $\Sig_{1}$-formula saying that `there exists an embedding
$f:X\to ON$', i.e.,
$\exi f\in{}^{\ome}ON\fal \bet,\gam\in X(\bet<\gam \to f(\bet)<f(\gam))$.

\eenu
}
\end{definition}
Note that for $\alpha\in OT(\Pi^{1}_{1})$,
$W(X)\cap\alpha=W(X\cap\alpha)$.
Also ${\sf KP}\ome\vdash \alp\in W_{\Sig}(X)\Rarw \alp\in W(X)$, and
${\sf KP}\ell\vdash \alp\in W(X) \Rarw \alp\in W_{\Sig}(X)$.

\begin{definition} \label{df:3wfdtg32}
{\rm 
For $X\subset OT(\Pi_{1}^{1})$ and
$\alpha\in OT(\Pi_{1}^{1})$,
\benu
\item
$
D[X] :\Lrarw 
\forall\alpha(\alpha\leq X\to W(\mathcal{C}^{\alpha}(X))\cap\alpha^{+}= X\cap\alpha^{+}
)$.

A set $X$ is said to be a \textit{distinguished set} if $D[X]$. 

\item\label{df:3wfdtg.832}
$D_{\Sig}[X]$ is a $\Sig$-formula defined by
\begin{equation}\label{eq:distinguishedclass}
D_{\Sig}[X] :\Lrarw 
\forall\alpha(\alpha\leq X\to W(\mathcal{C}^{\alpha}(X))\cap\alpha^{+}\subset X\cap\alpha^{+} \subset
W_{\Sig}(\mathcal{C}^{\alpha}(X))\cap\alpha^{+}
)
\end{equation}


\item\label{df:3wfdtg.932}
$\mathcal{W}:=\bigcup\{X: D_{\Sig}[X]\}$.


\eenu
}
\end{definition}

From ${\sf KP}\ome\vdash \alp\in W_{\Sig}(X)\Rarw \alp\in W(X)$ we see
$D_{\Sig}[X]\Rarw D[X]$ for any $X$.

Let $\alpha\in X$ for a $\Sig$-distinguished set $X$. 
Then $W(\mathcal{C}^{\alpha}(X))\cap\alpha^{+}= X\cap\alpha^{+}$.
Hence $X$ is a well order.
Although $\bigcup\{X: D[X]\}$ might be a proper class,
it turns out that $\calw$ is a set.


\bprp\label{prp:absolute}
Let $X\in L_{\mS}$.
\benu
\item\label{prp:absolute.0}
$\alp\in W(X)\Lrarw L_{\mS}\models \alp\in W(X)$.

\item\label{prp:absolute.1}
$\alp\in W_{\Sig}(X)\Lrarw L_{\mS}\models \alp\in W_{\Sig}(X)$.

\item\label{prp:absolute.2}
$D_{\Sig}[X]\Lrarw L_{\mS}\models D_{\Sig}[X]$.

\item\label{prp:absolute.3}
$D_{\Sig}[X] \lrarw D[X]$.

\item\label{prp:absolute.4}
$\mathcal{W}=\bigcup\{X\in L_{\mS} :D[X]\}$ and $\exi X(X=\calw)$.
\eenu
\eprp
\bprf
\ref{prp:absolute}.\ref{prp:absolute.0}.
Since $\alp\in W(X)$ is a $\Pi_{1}$-formula, it suffices to show 
$\alp\in W(X)$ assuming $L_{\mS}\models \alp\in W(X)$.
We obtain
$L_{\mS}\models \left(\alp\in W(X) \lrarw \alp\in W_{\Sig}(X)\right)$ by $L_{\mS}\models{\sf KP}\ell$.
Hence $\alp\in W_{\Sig}(X)$ and $\alp\in W(X)$.
\\
\ref{prp:absolute}.\ref{prp:absolute.1}.
Assume $\alp\in W_{\Sig}(X)$. 
Since $\alp\in W_{\Sig}(X)$ is a $\Sig_{1}$-formula, we obtain 
$L_{\mS}\models \alp\in W_{\Sig}(X)$ by $L_{\mS}\prec_{\Sig_{1}}L$.
The other direction follows from the persistency of $\Sig_{1}$-formulas.
\\
\ref{prp:absolute}.\ref{prp:absolute.2} follows from 
Propositions \ref{prp:absolute}.\ref{prp:absolute.0}
and \ref{prp:absolute}.\ref{prp:absolute.1}.
\\
\ref{prp:absolute}.\ref{prp:absolute.3}.
From $W_{\Sig}(\mathcal{C}^{\alpha}(X))\subset W(\mathcal{C}^{\alpha}(X))$
we see 
$D_{\Sig}[X] \to D[X]$.
Assume $D[X]$, $\alpha\leq X$ and $\bet\in X\cap\alp^{+}$.
Then $\bet\in W(\mathcal{C}^{\alpha}(X))\cap\alpha^{+}$ by $D[X]$.
We obtain $\bet\in W_{\Sig}(\mathcal{C}^{\alpha}(X))\cap\alp^{+}$
by $L_{\mS}\models \bet\in W(\mathcal{C}^{\alpha}(X))\to \bet\in W_{\Sig}(\mathcal{C}^{\alpha}(X))$ and Propositions \ref{prp:absolute}.\ref{prp:absolute.0}
and \ref{prp:absolute}.\ref{prp:absolute.1}.
\\
\ref{prp:absolute}.\ref{prp:absolute.4}.
By Proposition \ref{prp:absolute}.\ref{prp:absolute.3} we obtain
$\bigcup\{X\in L_{\mS} : D[X]\}\subset\calw$.
Let $\alp\in\calw$. This means a $\Sig_{1}$-formula
$\exi X(\alp\in X\land D_{\Sig}[X])$ holds. We obtain $L_{\mS}\models \exi X(\alp\in X\land D_{\Sig}[X])$ by 
$L_{\mS}\prec_{\Sig_{1}}L$.
By Propositions \ref{prp:absolute}.\ref{prp:absolute.2} and \ref{prp:absolute}.\ref{prp:absolute.3}
we obtain
$\alp\in\bigcup\{X\in L_{\mS} : L_{\mS}\models D_{\Sig}[X]\}=\bigcup\{X\in L_{\mS} : D_{\Sig}[X]\}=
\bigcup\{X\in L_{\mS} : D[X]\}$.
$\Del_{0}$-separation yields $\exi X(X=\calw)$.
\eprf

\begin{proposition}\label{lem:3.11.632}
Let $X\in L_{\mS}$ be a distinguished set. Then
$\alpha\in X \Rarw \forall\beta[\alpha\in\mathcal{C}^{\beta}(X)]$.
\end{proposition}

\begin{proposition}\label{lem:5uv.232general}
For any distinguished sets $X$ and $Y$ in $L_{\mS}$, 
$
X\cap\alpha=Y\cap\alpha \Rarw \forall\beta<\alpha^{+}\{\mathcal{C}^{\beta}(X)\cap\beta^{+}=\mathcal{C}^{\beta}(Y)\cap\beta^{+}\}$ holds
\end{proposition}

\begin{proposition}\label{lem:3wf5}
For distinguished sets $X$ and $Y$ in $L_{\mS}$,
$\alpha\leq X \spand \alpha\leq Y \Rarw X\cap\alpha^+=Y\cap\alpha^+$.
\end{proposition}

\begin{proposition}\label{lem:3wf6}
$\mathcal{W}$ is the maximal distinguished {\rm set}, i.e.,
$D[\mathcal{W}]$ and $\exists X(X=\mathcal{W})$.
\end{proposition}
\bprf
First we show $\fal\gam\in\calw(\gam\in\mathcal{C}^{\gam}(\calw))$.
Let $\gam\in\calw$, and pick a distinguished set $X\in L_{\mS}$ such that $\gam\in X$ by Proposition \ref{prp:absolute}.\ref{prp:absolute.4}.
Then $\gam\in\mathcal{C}^{\gam}(X)\subset\mathcal{C}^{\gam}(\calw)$ by $X\subset\calw$.

Let $\alp\leq\calw$. Pick a distinguished set $X\in L_{\mS}$ such that $\alp\leq X$.
We claim that $\calw\cap\alp^{+}=X\cap\alp^{+}$.
Let $Y\in L_{\mS}$ be a distinguished set and $\bet\in Y\cap\alp^{+}$.
Then $\bet\in Y\cap\bet^{+}=X\cap\bet^{+}$ by Proposition \ref{lem:3wf5}.
The claim yields
$W(\mathcal{C}^{\alp}(\calw))\cap\alp^{+}=W(\mathcal{C}^{\alp}(X))\cap\alp^{+}=
X\cap\alp^{+}=\calw\cap\alp^{+}$.
Hence $D[\mathcal{W}]$.
\eprf
\\

From $\calw\cap\alp^{+}=X\cap\alp^{+}$ for a $\Sig$-distinguished set $X$,
We see $\calw\cap\alp^{+}=X\cap\alp^{+}\subset W_{\Sig}(\mathcal{C}^{\alp}(X))\cap\alp^{+}=W_{\Sig}(\mathcal{C}^{\alp}(\calw))\cap\alp^{+}$.
Hence
$D_{\Sig}[\mathcal{W}]$.

\subsection{Sets $\mathcal{C}^{\alpha}(\mathcal{W}_{\mS})$ and $\mathcal{G}$}\label{subsec:C(X)}



\begin{definition}\label{df:calg}
$\mathcal{G}(Y):=\{\alpha\in OT(\Pi_{1}^{1}) :\alpha\in \mathcal{C}^{\alpha}(Y) \spand \mathcal{C}^{\alpha}(Y)\cap\alpha\subset Y\}$.
\end{definition}

\begin{lemma}\label{lem:wf5.332}
For $D[X]$,
$X\subset\mathcal{G}(X)$.
\end{lemma}

\blem\label{lem:6.21}
Suppose $D[Y]$ and $\alp\in\mathcal{G}(Y)$ for $Y\in L_{\mS}$.
Let 
$X=W(\mathcal{C}^{\alp}(Y))\cap\alp^{+}\in L_{\mS}$.
Assume that the following condition (\ref{eq:6.21.55})
is fulfilled.
Then $\alp\in X$ and $D[X]$.

\beqn
\fal\bet<\mS\left(
Y\cap\alp^{+}<\bet \spand \bet^{+}<\alp^{+} \to
W(\mathcal{C}^{\bet}(Y))\cap\bet^{+}\subset Y
\right)
\label{eq:6.21.55}
\eeqn

\elem

\bprp\label{prp:noncritical}
Let $D[X]$.
\benu
\item\label{prp:noncritical.1}
Let $\{\alp,\bet\}\subset X$ with $\alp+\bet=\alp\#\bet$ and $\alp>0$.
Then $\gam=\alp+\bet\in X$.

\item\label{prp:noncritical.2}
If $\{\alp,\bet\}\subset X$, then $\vphi_{\alp}(\bet)\in X$.

\eenu
\eprp

\bprp\label{lem:6.29}
\benu
\item\label{lem:6.29.1}
$0\in\calw$.

\item\label{lem:6.29.2}
Let either $\sig=0$ or $\sig=\psi_{\mS}^{f}(a)$ or $\sig=\psi_{\pi}^{f}(a)$.
Assume $\sig\in \calw$. Then $\sig^{+}\in\calw$.
\eenu
\eprp
\bprf
Each is seen from Lemma \ref{lem:6.21} as follows.
\\
\ref{lem:6.29}.\ref{lem:6.29.1}.
We see $0\in Y=W(\calc^{0}(\emptyset))\cap \Ome\in L_{\mS}$ with $\Ome=0^{+}$ and $D[Y]$.
\\
\ref{lem:6.29}.\ref{lem:6.29.2}.
Let $\sig\in Y\in L_{\mS}$ with $D[Y]$.
We see $\sig^{+}\in X=W(\calc^{\sig^{+}}(Y))\cap\sig^{++}\in L_{\mS}$ and $D[X]$.
\eprf

\begin{lemma}\label{th:3wf16}
Suppose $D[Y]$ 
with $\{0,\Ome\}\subset Y\in L_{\mS}$, and for $\eta\in OT(\Pi_{1}^{1})\cap(\mS+1)$
\begin{equation}\label{eq:3wf16hyp.132}
\eta\in\mathcal{G}(Y)
\end{equation}
and
\begin{equation}\label{eq:3wf16hyp.232}
\fal\gam\prec\eta (\gam\in\mathcal{G}(Y)
\Rarw \gam\in Y)
\end{equation}
Let $X=W(\mathcal{C}^{\eta}(Y))\cap\eta^{+}$.
Then
$\eta\in X\in L_{\mS}$ and $D[X]$.
\end{lemma}

\subsection{Mahlo universes}
In this subsection we consider the maximal distinguished class $\calw^{P}$ \textit{inside} a set
$P\in L_{\mS}$ as in \cite{Wienpi3d}.
Let $ad$ denote a $\Pi^{-}_{3}$-sentence such that a transitive set $z$ is admissible 
iff $(z;\in)\models ad$.
Let $lmtad :\Lrarw \fal x\exi y(x\in y \land ad^{y})$.
Observe that $lmtad$ is a $\Pi^{-}_{2}$-sentence.

\bdf\label{df:3auni}
{\rm
\benu
\item 
By a \textit{universe} 
we mean either $L_{\mS}$ or a transitive set $Q\in L_{\mS}$ with $\ome\in Q$. 
Universes are denoted by $P,Q,\ldots$

\item For a universe $P$ and a set-theoretic sentence $\vphi$, 
$P\models\vphi :\Lrarw (P;\in)\models\vphi$.

\item 
A universe $P$ is said to be a \textit{limit universe} if $lmtad^{P}$ holds, i.e., $P$
 is a limit of admissible sets. The class of limit universes is denoted by $Lmtad$.
 
\eenu
}
\edf

\blem\label{lem:3ahier} 
$W(\mathcal{C}^{\alp}(X))$ as well as $D[X]$ are absolute for limit universes $P$.


\elem

\bdf\label{df:3awp}
{\rm 
For a universe $P$, let
$
\calw^{P}:=\bigcup\{X\in P:D[X]\}
$.
}
\edf

We see $\calw=\calw^{L_{\mS}}$ from Proposition \ref{prp:absolute}.\ref{prp:absolute.4}.

\blem\label{lem:3afin}
Let $P$ be a universe closed under finite unions, and $\alp\in OT(\Pi_{1}^{1})$.
\benu
\item\label{lem:3afin.1}
There is a finite set $K(\alp)\subset OT(\Pi_{1}^{1})$ such that
$
\fal Y\in P\fal\gam
[
K(\alp)\cap Y=K(\alp)\cap \calw^{P} 
\Rarw 
\left(
\alp\in \mathcal{C}^{\gam}(\calw^{P})
\Lrarw \alp\in \mathcal{C}^{\gam}(Y)
\right)
]
$.

\item\label{lem:3afin.2}
There exists a distinguished set $X\in P$ such that
$
\fal Y\in P\fal\gam
[
X\subset Y\spand D[Y]
\Rarw (\alp\in \mathcal{C}^{\gam}(\calw^{P})
\Lrarw \alp\in \mathcal{C}^{\gam}(Y)
)
]
$.
\eenu
\elem

\begin{proposition}\label{lem:3wf6}
For each limit universe $P$,
$D[\mathcal{W}^{P}]$ holds, and
 $\exists X\in L_{\mS}(X=\mathcal{W}^{P})$ if $P\in L_{\mS}$.
\end{proposition}

For a universal $\Pi_{n}$-formula $\Pi_{n}(a)\, (n>0)$ uniformly on admissibles, let
\[
P\in M_{2}(\calc) :\Lrarw P\in Lmtad \spand 
\fal b\in P[P\models\Pi_{2}(b)\rarw \exi Q\in \calc\cap P(Q\models\Pi_{2}(b))]
.\]

\bdf\label{df:bigO}
{\rm
Let $\gam=\psi_{\sigma}^{f}(\alpha)\prec\mS$ and
$\mu=\Lam(\gam)<\mS^{+}$ be the ordinal in Definition \ref{df:Lammu}.
Let $O(\gam)=o_{\mu}(f)<\mS^{+}$, where 
$o_{\mu}(f)$ is the ordinal defined in Definition \ref{df:Lam.of} from the epsilon number $\mu$.

Let
$O(\Ome)=1$, $O(\mS)=\mS^{+}$ and $O(\gam)=0$ else.
}
\edf

\newpage

\blem\label{lem:4acalg} 
Let $\calc$ be a $\Pi^{1}_{0}$-class such that $\calc\subset Lmtad$.
Suppose $P\in M_{2}(\calc)$, $\alp\in\calg(\calw^{P})$ and $O(\alp)\in W(\mathcal{C}^{\mS}(\calw^{P}))$
Then there exists a universe $Q\in \calc$ such that $\alp\in\calg(\calw^{Q})$
and $O(\alp)\in W(\mathcal{C}^{\mS}(\calw^{Q}))$. 
\elem
\bprf 
Suppose $P\in M_{2}(\calc)$, $\alp\in\calg(\calw^{P})$ and $O(\alp)\in W(\mathcal{C}^{\mS}(\calw^{P}))$.
First by $\alp\in \mathcal{C}^{\alp}(\calw^{P})$,
$O(\alp)\in\mathcal{C}^{\mS}(\calw^{P})$ and Lemma \ref{lem:3afin} pick a distinguished set
 $X_{0}\in P$ such that $\alp\in \mathcal{C}^{\alp}(X_{0})$,
 $O(\alp)\in\mathcal{C}^{\mS}(X_{0})$ and
$K(\alp)\cap\calw^{P}\subset X_{0}$.
Then for any universe $X_{0}\in Q\in P$, we obtain $O(\alp)\in W(\mathcal{C}^{\mS}(\calw^{Q}))$
by $\calw^{Q}\subset\calw^{P}$.

Next writing $\mathcal{C}^{\alp}(\calw^{P})\cap\alp\subset\calw^{P}$ analytically we have

\[
\fal\bet<\alp[\bet\in \mathcal{C}^{\alp}(\calw^{P}) \Rarw 
\exi Y\in P(D[Y]\spand \bet\in Y)]
\]
By Lemma \ref{lem:3afin} we obtain
$\bet\in \mathcal{C}^{\alp}(\calw^{P})\Lrarw
\exi X\in P\{D[X] \spand K(\bet)\cap\calw^{P}\subset X \spand \bet\in \mathcal{C}^{\alp}(X) \}$.
Hence
for any $\bet<\alp$ and any distinguished set $X\in P$,
 there are $\gam\in K(\bet)$, $Z\in P$ and a distinguished set $Y\in P$ such that 
 if $\gam\in Z \spand D[Z] \to \gam\in X$ and $\bet\in \mathcal{C}^{\alp}(X)$, then 
 $\bet\in Y$.
By Lemma \ref{lem:3ahier} $D[X]$ is absolute for limit universes. 
Hence the following $\Pi_{2}$-predicate holds in the universe $P\in M_{2}(\calc)$:
\beqnarr
&& \fal\bet<\alp\fal X\exi\gam\in K(\bet)\exi Z\exi Y
[
\{
D[X]\spand 
 (\gam\in Z \spand D[Z] \to \gam\in X)
  \spand \bet\in \mathcal{C}^{\alp}(X)
\}
\nonumber
\\
&&
 \Rarw 
 \left(
 D[Y]\spand \bet\in Y
 \right)
 ]
 \label{eq:3acalg}
\eeqnarr

Now pick a universe $Q\in \calc\cap P$ with $X_{0}\in Q$ and $Q\models (\ref{eq:3acalg})$. 
Tracing the above argument backwards in the limit universe $Q$ we obtain 
$\mathcal{C}^{\alp}(\calw^{Q})\cap\alp\subset\calw^{Q}$ and 
$X_{0}\subset\calw^{Q}=\bigcup\{X\in Q: Q\models D[X]\}\in P$. 
Thus Lemma \ref{lem:3afin} yields $\alp\in \mathcal{C}^{\alp}(\calw^{Q})$. 
We obtain $\alp\in\calg(\calw^{Q})$.
\eprf

\begin{proposition}\label{prp:psimS}
Let $\gam=\psi_{\sigma}^{f}(\alpha)\in \mathcal{G}(Y)$ and
$\gam\preceq\gam_{0}=\psi_{\mS}^{g}(b)$ with $b={\tt p}_{0}(\gam)$.
Then $O(\gam)\in \calc^{\mS}(Y)$.
\end{proposition}
\bprf
We have
$\gamma\in\mathcal{C}^{\gamma}(Y)$ and 
$\mathcal{C}^{\gamma}(Y)\cap\gamma\subset Y$.
We obtain $\{\sig,b\}\cup SC_{\mS^{+}}(f)\subset\mathcal{C}^{\gamma}(Y)$.
We obtain $SC_{\mS^{+}}(f,b)<\mu=\Lam(\gam)$ by Definition \ref{df:notationsystempi11}. 
$E_{\mS}(SC_{\mS^{+}}(f,b))\subset\mathcal{C}^{\gamma}(Y)$ follows from $\gam<\mS$.
On the other hand we have $SC_{\mS^{+}}(f,b)\subset\calh_{b}(\gam)$ for $b={\tt p}_{0}(\gam)$
by (\ref{eq:notationsystem.11}).
This yields $E_{\mS}(SC_{\mS^{+}}(f,b))\subset\calh_{b}(\gam)\cap\mS\subset\gam$.
We obtain
$E_{\mathbb{S}}(SC_{\mS^{+}}(f,b))\subset\mathcal{C}^{\gamma}(Y)\cap\gamma\subset Y$.
Hence
$SC_{\mS^{+}}(f,b)\subset\mathcal{C}^{\mathbb{S}}(Y)$.
From $SC_{\mS^{+}}(f,b)\subset \calc^{\mS}(Y)$ we see
 $O(\gam)=o_{\mu}(f)\in \calc^{\mS}(Y)$ for $\gamma=\psi_{\sigma}^{f}(\alpha)\in\mathcal{G}(Y)$.
\eprf

\bdf\label{df:UVM}
{\rm 
We define the class $M_{2}(\alp)$ of $\alp$-recursively Mahlo universes for 
$\mS\geq\alp\in OT(\Pi_{1}^{1})$ as follows:
\beqn\label{eqarr:4a(1)}
 P\in M_{2}(\alp) \Lrarw 
 P\in Lmtad \spand 
 \fal\bet\prec\alp
 [
O(\bet)\in W(\calc^{\mS}(\calw^{P}))
 \Rarw P\in M_{2}(M_{2}(\bet))] 
\eeqn
$M_{2}(\alp)$ is a $\Pi_{3}$-class.

}
\edf

\blem\label{lem:3awf16.1} 
If $\mS\geq\eta\in\calg(\calw^{P})$, $O(\eta)\in W(\mathcal{C}^{\mS}(\calw^{P}))$
and $P\in M_{2}(M_{2}(\eta))$ with $P\in L_{\mS}$, 
then $\eta\in\calw^{P}$.
\elem
\bprf 
We show this by induction on $\in$. Suppose, as IH, the lemma holds for any $Q\in P$.
By Lemma \ref{lem:4acalg} pick a $Q\in P$ such that $Q\in M_{2}(\eta)$, and for $Y=\calw^{Q}\in P$,
$\{0,\Ome\}\subset Y$, $O(\eta)\in W(\mathcal{C}^{\mS}(Y))$ and
\beqn\renewcommand{\theequation}{\ref{eq:3wf16hyp.132}} 
\eta\in\calg(Y)
\eeqn
\addtocounter{equation}{-1}
On the other the definition (\ref{eqarr:4a(1)}) yields 
$\fal\gam\prec\eta[
O(\gam)\in W(\calc^{\mS}(\calw^{Q}))
 \Rarw Q\in M_{2}(M_{2}(\gam))]$.
IH yields with $Y=\calw^{Q}$
\[
\fal\gam\prec\eta(\gam\in\calg(Y)\spand O(\gam)\in W(\calc^{\mS}(Y)) \Rarw  \gam\in Y)
\]
On the other $O(\eta)\in W(\mathcal{C}^{\mS}(Y))$ yields $O(\gam)\in W(\calc^{\mS}(Y))$
for $\calg(Y)\ni\gam\prec\eta$ by Proposition \ref{prp:psimS}. Therefore
\beqn\renewcommand{\theequation}{\ref{eq:3wf16hyp.232}}
\fal\gam\prec\eta(\gam\in\calg(Y) \Rarw  \gam\in Y)
\eeqn
\addtocounter{equation}{-1}
Therefore by Lemma \ref{th:3wf16} we conclude 
$\eta\in X$ and $D[X]$ for $X=W(\mathcal{C}^{\eta}(Y))\cap\eta^{+}$.

$X\in P$ follows from $Y\in P\in Lmtad$. 
Consequently $\eta\in\calw^{P}$.
\eprf

\blem\label{lem:4aro}
 $\fal\eta\leq\mS[O(\eta)\in W(\mathcal{C}^{\mS}(\calw)) \Rarw L_{\mS}\in M_{2}(M_{2}(\eta))]$.
\elem
\bprf 
We show the lemma by induction on $O(\eta)\in W(\mathcal{C}^{\mS}(\calw))$.
Suppose $O(\eta)\in W(\mathcal{C}^{\mS}(\calw))$ and $L_{\mS}\models\Pi_{2}(b)$ for a $b\in L_{\mS}$. 
We have to find a universe $Q\in L_{\mS}$ such that $b\in Q$, $Q\in M_{2}(\eta)$ and $Q\models\Pi_{2}(b)$.

By the definition (\ref{eqarr:4a(1)}) $L_{\mS}\in M_{2}(\eta)$ is equivalent to
$\fal\gam\prec\eta [ O(\gam)\in W(\mathcal{C}^{\mS}(\calw)) \Rarw L_{\mS}\in M_{2}(M_{2}(\gam))]$, where
$\calw=\calw^{L_{\mS}}$ by Proposition \ref{prp:absolute}.\ref{prp:absolute.4}.
We obtain $\gam\prec\eta \Rarw O(\gam)<O(\eta)$. 
Thus IH yields $L_{\mS}\in M_{2}(\eta)$.
Let $g$ be a primitive recursive function in the sense of set theory such that 
$L\in M_{2}(\eta) \Lrarw P\models\Pi_{3}(g(\eta))$.
Then $L_{\mS}\models\Pi_{2}(b) \land \Pi_{3}(g(\eta))$. 
Since this is a $\Pi_{3}$-formula which holds in a $\Pi_{3}$-reflecting universe $L_{\mS}$, 
we conclude for some $Q\in L_{\mS}$, $Q\models\Pi_{2}(b)\land \Pi_{3}(g(\eta))$
 and hence $Q\in M_{2}(\eta)$. We are done.
\eprf

\blem\label{lem:3awf16.2}
$\fal\eta\leq\mS\left[
\eta\in\calg(\calw) \spand O(\eta)\in W(\mathcal{C}^{\mS}(\calw)) \Rarw \eta\in\calw\right]$.
\elem
\bprf 
Assume $\mS\geq\eta\in\calg(\calw)$ and $O(\eta)\in W(\mathcal{C}^{\mS}(\calw))$.
Lemma \ref{lem:4aro} yields $L_{\mS}\in M_{2}(M_{2}(\eta))$. 
From this we see $L_{\mS}\in M_{2}(\calc)$ with $\calc=M_{2}(M_{2}(\eta))$ 
as in the proof of Lemma \ref{lem:4aro} using $\Pi_{3}$-reflection of $L_{\mS}$ once again. 
Then by Lemma \ref{lem:4acalg} pick a set $P\in L_{\mS}$ such that 
$\eta\in\calg(\calw^{P})$, $O(\eta)\in W(\mathcal{C}^{\mS}(\calw^{P}))$ and $P\in\calc=M_{2}(M_{2}(\eta))$. 
Lemma \ref{lem:3awf16.1} yields $\eta\in\calw^{P}\subset\calw$.
\eprf

\bdf
{\rm Let
$\calw_{1}:=W(\calc^{\mS}(\calw))$.
}
\edf

\bprp\label{prp:w1}
\benu
\item\label{prp:w1.1}
$\calc^{\mS^{+}}(\calw_{1})\cap\mS^{+}=\calw_{1}\cap\mS^{+}$ and
$\calc^{\mS^{+}}(\calw_{1})\cap\mS=\mathcal{C}^{\mathbb{S}}(\mathcal{W})\cap\mathbb{S}=\calw_{1}\cap\mS=\calw\cap\mS$.

\item\label{prp:w1.2}
$\mS\in \calw_{1}$.

\item\label{prp:w1.3}
$TI[\calc^{\mS^{+}}(\calw_{1})\cap\ome_{n}(\mS^{+}+1)]$ for {\rm each} $n<\ome$.
\eenu
\eprp
\bprf
\ref{prp:w1}.\ref{prp:w1.1} and \ref{prp:w1}.\ref{prp:w1.2}.
Since there is no regular ordinal$>\mS^{+}$,
$\calc^{\mS^{+}}(\calw_{1})\cap\mS^{+}=\calw_{1}\cap\mS^{+}$.
We see $\mathcal{C}^{\mathbb{S}}(\mathcal{W})\cap\mathbb{S}=\mathcal{W}\cap\mS=\calw_{1}\cap\mS$ 
from $\psi_{\mS^{+}}(a)>\mS$ and $W(\calw)=\calw$.
Hence $\mS\in \calw_{1}$.
\\
\ref{prp:w1}.\ref{prp:w1.3}.
$TI[\calc^{\mS^{+}}(\calw_{1})\cap\mS^{+}]$ follows from
$\calc^{\mS^{+}}(\calw_{1})\cap\mS^{+}=\calw_{1}\cap\mS^{+}$.
By meta-induction on $n<\ome$, we see $TI[\calc^{\mS^{+}}(\calw_{1})\cap\ome_{n}(\mS^{+}+1)]$
using the Gentzen's jump set.
\eprf

\blem\label{th:id5wf21S+}
$a\in\calc^{\mS^{+}}(\calw_{1})\cap\ome_{n}(\mS^{+}+1)\spand
\psi_{\mS^{+}}(a)\in OT(\Pi^{1}_{1}) \Rarw \psi_{\mS^{+}}(a)\in\calw_{1}$
for {\rm each} $n<\ome$.
\elem
\bprf
By Proposition \ref{prp:w1}.\ref{prp:w1.3} it suffices to show that
$Prg[\calc^{\mS^{+}}(\calw_{1}), B]$ for
$B(a):\Lrarw \left[\psi_{\mS^{+}}(a)\in OT(\Pi^{1}_{1}) \Rarw \psi_{\mS^{+}}(a)\in\calw_{1}\right]$.
Assume $a\in \calc^{\mS^{+}}(\calw_{1})$.
We obtain $\psi_{\mS^{+}}(a)\in \calc^{\mS}(\calw)$ 
by Propositions \ref{lem:CX2}.\ref{lem:CX2.3} and \ref{prp:w1}.\ref{prp:w1.1}.

Next we show 
$\bet\in \calc^{\mS}(\calw)\cap\calh_{a}(\psi_{\mS^{+}}(a)) \Rarw \bet\in\calc^{\mS^{+}}(\calw_{1})$
by induction on $\ell\bet$.
By Proposition \ref{prp:w1} we may assume $\bet=\psi_{\mS^{+}}(b)$.
Then $b\in \calc^{\mS}(\calw)\cap\calh_{a}(\psi_{\mS^{+}}(a))$, and $b\in \calc^{\mS^{+}}(\calw_{1})$
by IH on lengths.
Moreover $b<a$. Hence IH yields $\bet\in\calw_{1}\cap\mS^{+}\subset\calc^{\mS^{+}}(\calw_{1})$.

In particular we obtain
$\calc^{\mS}(\calw)\cap\psi_{\mS^{+}}(a)\subset\calc^{\mS^{+}}(\calw_{1})\cap\mS^{+}\subset\calw_{1}$.
Therefore $\psi_{\mS^{+}}(a)\in W(\calc^{\mS}(\calw))=\calw_{1}$.
\eprf

\subsection{Well-foundedness proof concluded}

\begin{definition}\label{df:id4wfA}
{\rm 
For irreducible functions $f$ let
\[
f\in J:\Lrarw SC_{\mS^{+}}(f)\subset \calw_{1}.
\]

For $a\in OT(\Pi^{1}_{1})$ and irreducible functions $f$, define:
\beqnarrs
 A(a,f) & :\Lrarw &
 \forall\sigma\in \calw_{1}\cap\mS^{+}[\psi_{\sigma}^{f}(a)\in OT(\Pi^{1}_{1}) 
 \spand O(\psi_{\sigma}^{f}(a))\in\calw_{1} \Rarw \psi_{\sigma}^{f}(a)\in\mathcal{W}].
\\
\mbox{{\rm MIH}}(a) & :\Lrarw &
 \forall b\in \calc^{\mS^{+}}(\calw_{1})\cap a\forall f\in J \, A(b,f).
\\
\mbox{{\rm SIH}}(a,f) & :\Lrarw &
 \forall g\in J [g<^{0}_{lx}f  \Rarw A(a,g)].
\eeqnarrs
}
\end{definition}

\begin{lemma}\label{th:id5wf21}
Assume $a\in \calc^{\mS^{+}}(\calw_{1})\cap\ome_{n}(\mS^{+}+1)$, $f\in J$, $\mbox{{\rm MIH}}(a)$, and $\mbox{{\rm SIH}}(a,f)$ in Definition \ref{df:id4wfA}.
Then
\[
 \forall\kappa\in\mathcal{W}_{1}\cap\mS^{+}[\psi_{\kappa}^{f}(a)\in OT(\Pi_{1}^{1})\spand
 O(\psi_{\kappa}^{f}(a))\in\calw_{1}  \Rarw \psi_{\kappa}^{f}(a)\in\mathcal{W}].
\]
\end{lemma}
\bprf
This is seen as in \cite{Wienpi3d,singlewfprf} from
Lemma \ref{lem:3awf16.2}.
Let $\alpha_{1}=\psi_{\kappa}^{f}(a)\in OT(\Pi_{1}^{1})$ be such that
$O(\alp_{1})\in\calw_{1}$,
 $a\in\calc^{\mS^{+}}(\calw_{1})$,
$\mS\geq\kappa\in\mathcal{W}_{1}$ and $f\in J$. 
By Lemma \ref{lem:3awf16.2} it suffices to show $\alpha_{1}\in\mathcal{G}(\calw)$.

By Proposition \ref{lem:CX2}.\ref{lem:CX2.3} we have
$\{\kappa,a\}\cup SC_{\mS^{+}}(f)\subset\mathcal{C}^{\alpha_{1}}(\mathcal{W})$, and hence $\alpha_{1}\in\mathcal{C}^{\alpha_{1}}(\mathcal{W})$.
It suffices to show the following claim.
\begin{equation}\label{clm:id5wf21.1}
\forall\beta_{1}\in\mathcal{C}^{\alpha_{1}}(\mathcal{W})\cap\alpha_{1}[\beta_1\in\mathcal{W}].
\end{equation}
\textbf{Proof} of (\ref{clm:id5wf21.1}) by induction on $\ell\beta_1$. 
Assume $\beta_{1}\in\mathcal{C}^{\alpha_{1}}(\mathcal{W})\cap\alpha_{1}$ and let
\[
\mbox{LIH} :\Lrarw
\forall\gamma\in\mathcal{C}^{\alpha_{1}}(\mathcal{W})\cap\alpha_{1}[\ell\gamma<\ell\beta_{1} \Rarw \gamma\in\mathcal{W}].
\]

We show $\beta_1\in\mathcal{W}$. 
We may assume that $\beta_{1}=\psi_{\pi}^{g}(b)$ for some $\pi,b,g$
such that $\{\pi,b\}\cup SC_{\mS^{+}}(g)\subset\mathcal{C}^{\alpha_{1}}(\mathcal{W})$ and
$\alp_{1}<\pi\leq\mS$.
\\
\textbf{Case 1}.
$b<a$, $\beta_{1}<\kappa$ and $\{\pi,b\}\cup SC_{\mS^{+}}(g)\subset\calh_{a}(\alp_{1})$:
Let $B$ denote a set of subterms of $\beta_{1}$ defined recursively as follows.
First $\{\pi,b\}\cup SC_{\mS^{+}}(g)\subset B$.
Let $\alpha_{1}\leq\beta\in B$. 
If $\beta=_{NF}\gamma_{m}+\cdots+\gamma_{0}$, then $\{\gamma_{i}:i\leq m\}\subset B$.
If $\beta=_{NF}\varphi\gamma\delta$, then $\{\gamma,\delta\}\subset B$.
If $\beta=_{NF}\gam^{+}$, then $\gamma\in B$.
If $\beta=\psi_{\sigma}^{h}(c)$ with $\sig>\alp_{1}$, then $\{\sigma,c\}\cup SC_{\mS^{+}}(h)\subset B$.

Then from $\{\pi,b\}\cup SC_{\mS^{+}}(g)\subset\mathcal{C}^{\alpha_{1}}(\mathcal{W})$ we see inductively that
$B\subset\mathcal{C}^{\alpha_{1}}(\mathcal{W})$.
Hence by LIH we obtain $B\cap\alpha_{1}\subset\mathcal{W}$.
Moreover if $\alpha_{1}\leq\psi_{\sigma}^{h}(c)\in B$, 
then $c<a$.
We claim that
\begin{equation}\label{eq:case2A}
\forall\beta\in B(\beta\in\calc^{\mS^{+}}(\mathcal{W}_{1}))
\end{equation}
\textbf{Proof} of (\ref{eq:case2A}) by induction on $\ell\beta$.
Let $\beta\in B$. 
We can assume that $\alpha_{1}\leq\beta=\psi_{\sigma}^{h}(c)$ by LIH.
Then by induction hypothesis we have 
$\{\sigma,c\}\cup SC_{\mS^{+}}(h)\subset\calc^{\mS^{+}}(\mathcal{W}_{1})$.
On the other hand we have $c<a<\ome_{n}(\mS^{+}+1)$.
If $\sig=\mS^{+}$, then Lemma \ref{th:id5wf21S+} yields $\bet=\psi_{\mS^{+}}(c)\in\calw_{1}$.
Let $\sig\leq\mS$. Then 
$\{\sigma\}\cup SC_{\mS^{+}}(h)\subset\calc^{\mS^{+}}(\mathcal{W}_{1})\cap\mS^{+}=\calw_{1}\cap\mS^{+}$.
Let $\alp_{1}\leq\bet=\psi_{\sig}^{h}(c)\preceq\psi_{\mS}^{h_{0}}(c_{0})$.
Then $SC_{\mS^{+}}(c_{0})\subset B$ and $SC_{\mS^{+}}(c_{0})\subset\calw_{1}$ by
induction hypothesis on lengths.
Hence $\Lam(\bet)\in\calw_{1}$ for the least epsilon number $\Lam(\bet)>\max(SC_{\mS^{+}}(c_{0}))$.
We obtain $O(\bet)\in\calw_{1}$ by 
$SC_{\mS^{+}}(h)\cup\{\Lam(\bet)\}\subset\calw_{1}$.
$\mbox{MIH}(a)$ yields $\beta\in\mathcal{W}$.
Thus (\ref{eq:case2A}) is shown.
\hspace*{\fill} $\Box$
\\

In particular we obtain $\{\pi,b,\Lam(\bet_{1})\}\cup SC_{\mS^{+}}(g)\subset \calc^{\mS^{+}}(\mathcal{W}_{1})$.
Moreover we have $b<a$.
Therefore once again $\mbox{MIH}(a)$ yields $\beta_{1}\in\mathcal{W}$.
\\
\textbf{Case 2}. 
$b=a$, $\pi=\kappa$, $SC_{\mS^{+}}(g)\subset\calh_{a}(\alp_{1})$ and $g<^{0}_{lx}f$: 
As in (\ref{eq:case2A}) we see that $SC_{\mS^{+}}(g)\subset\mathcal{W}_{1}$ from 
Lemma \ref{th:id5wf21S+} and $\mbox{MIH}(a)$.
$\mbox{SIH}(a,f)$ yields $\beta_{1}\in\mathcal{W}$.
\\
\textbf{Case 3}.
$a\leq b$ and $SC_{\mS^{+}}(f)\cup\{\kappa,a\}\not\subset\calh_{b}(\bet_{1})$:
As in as in \cite{Wienpi3d,singlewfprf}  we see that
there exists a $\gamma$ such that $\beta_{1}\leq\gamma\in\mathcal{W}\cap\alpha_{1}$.
Then $\beta_{1}\in\mathcal{W}$ follows from $\beta_{1}\in\mathcal{C}^{\alpha_{1}}(\mathcal{W})$.

This completes a proof of (\ref{clm:id5wf21.1}) and of the lemma.
\hspace*{\fill} $\Box$

\begin{lemma}\label{lem:psiw}
For $\psi_{\kappa}^{f}(a)\in OT(\Pi^{1}_{1})$, if
$a\in \calc^{\mS^{+}}(\calw_{1})\cap\ome_{n}(\mS^{+}+1)$,
$\{\kappa\}\cup SC_{\mS^{+}}(f)\subset \calw_{1}\cap\mS^{+}$ and 
$O(\psi_{\kappa}^{f}(a))\in\calw_{1}$,
then $\psi_{\kappa}^{f}(a)\in \mathcal{W}$.
\end{lemma}
\bprf
This is seen from Lemma
\ref{th:id5wf21} and
Proposition \ref{prp:w1}.\ref{prp:w1.3}.
Note that if $\bet=\psi_{\kap}^{g}(a)<\psi_{\kap}^{f}(a)=\alp$ by $g<_{lx}^{0}f$, then
${\tt p}_{0}(\bet)={\tt p}_{0}(\alp)$, $\Lam(\bet)=\Lam(\alp)$ and $O(\bet)<O(\alp)$ by
Lemma \ref{lem:oflx}.
\hspace*{\fill} $\Box$

\blem\label{lem:eachpi11}
For {\rm each} $\alp\in OT(\Pi^{1}_{1})$,
$\alp\in\calc^{\mS^{+}}(\calw_{1})$.
\elem
\bprf
This is seen by meta-induction on $\ell\alp$ using 
Propositions \ref{prp:noncritical} and \ref{lem:6.29}, and
Lemmas \ref{th:id5wf21S+} and
\ref{lem:psiw}.
\eprf
\\

\noindent
\textbf{Proof} of Theorem \ref{th:wf}.
For each $\alp\in OT(\Pi^{1}_{1})$ we obtain $\alp\in\calc^{\mS^{+}}(\calw_{1})$ by Lemma \ref{lem:eachpi11}.
Therefore by Proposition \ref{prp:w1}.\ref{prp:w1.1} we obtain for each $n<\ome$,
$\psi_{\Ome}(\ome_{n}(\mS^{+}+1))\in \calc^{\mS^{+}}(\calw_{1})\cap\Ome=\calw\cap\Ome
=W(\calc^{0}(\emptyset))\cap\Ome$, where
$W(\calc^{0}(\emptyset))=W(OT(\Pi^{1}_{1}))$.


\begin{thebibliography}{99}

\bibitem{Wienpi3d} 
T. Arai, 
Wellfoundedness proofs by means of non-monotonic inductive definitions I: $\Pi^{0}_{2}$-operators, 
\textit{Jour. Symb. Logic} {\bf 69} (2004) 830--850.

\bibitem{singlewfprf}
T. Arai,
Wellfoundedness proof with the maximal distinguished set,
to appear in \textit{Arch. Math. Logic}.

\bibitem{singlestable}
T. Arai,
An ordinal analysis of a single stable ordinal,
submitted.

\bibitem{ghent}
T. Arai, Lectures on ordinal analysis,
a lecture notes for a mini-course in Department of Mathematics, Ghent University,
14 Mar.-25 Mar. 2023.


\bibitem{Buchholz75}
W. Buchholz, 
Normalfunktionen und konstruktive Systeme von Ordinalzahlen.
In: Diller, J., M\"uller, G. H. (eds.)
Proof Theory Symposion Kiel 1974, Lect. Notes Math. vol. 500, pp. 4-25, Springer (1975)


\bibitem{RathjenAFML2}M. Rathjen, 
An ordinal analysis of parameter free $\Pi^{1}_{2}$-comprehension,
\textit{Arch. Math. Logic} \textbf{44} (2005) 263-362.



\end{thebibliography}
\end{document}